\documentclass [12pt]{amsart}
\usepackage{amsmath,amssymb,amscd}
\usepackage{pstricks}

\def\e{\epsilon}

\def\bp{\begin{proposition}}
\def\ep{\end{proposition}}
\def\bt{\begin{theo}}
\def\et{\end{theo}}
\def\be{\begin{equation}}
\def\ee{\end{equation}}
\def\bl{\begin{lemma}}
\def\el{\end{lemma}}
\def\bc{\begin{corollary}}
\def\ec{\end{corollary}}
\def\pr{\noindent{\bf Proof: }}

\def\bd{\begin{definition}}
\def\ed{\end{definition}}

\def\dist{\text{dist}}
\begin{document}

\newtheorem{theo}{Theorem}[section]
\newtheorem{lemma}[theo]{Lemma}
\newtheorem{definition}[theo]{Definition}
\newtheorem{corollary}[theo]{Corollary}
\newtheorem{proposition}[theo]{Proposition}
\numberwithin{equation}{section}

\begin{titlepage}

\begin{center}

\topskip 5mm

{\LARGE{\bf {Linear versus Non-linear}}} \vskip 4mm {\LARGE{\bf
{Acquisition of Step-Functions}}} \vskip 8mm

{\large {\bf B. Ettinger$^{**}$, N. Sarig$^{*}$, Y. Yomdin $^{*}$}}

\vspace{6 mm}
\end{center}

{$^{*}$ Dept. of Math., The Weizmann Inst. of Science, Rehovot
76100, Israel. e-mail: niv.sarig@weizmann.ac.il,
yosef.yomdin@weizmann.ac.il}\\
{$^{**}$ Dept. of Math., Univ. of California, Berkeley, CA 94720.
e-mail: ettinger@math.berkeley.edu}

\vspace{4 mm}

\begin{center} {\it To G. Henkin on his 65-th birthday}

\end{center}

\vspace{4 mm}

\begin{center}

{ \bf Abstract}
\end{center}

{\small {We address in this paper the following two closely
related problems:

\noindent 1. How to represent functions with singularities (up to
a prescribed accuracy) in a compact way? 2. How to reconstruct
such functions from a small number of measurements? The stress is
on a comparison of linear and non-linear approaches. As a model
case we use piecewise-constant functions on $[0,1]$, in
particular, the Heaviside jump function
$\mathcal{H}_t=\chi_{[0,t]}.$ Considered as a curve in the Hilbert
space $L^2([0,1])$  it is completely characterized by the fact
that any two its disjoint chords are orthogonal. We reinterpret
this fact in a context of step-functions in one or two variables.

Next we study the limitations on representability and
reconstruction of piecewise-constant functions by linear and
semi-linear methods. Our main tools in this problem are
Kolmogorov's $n$-width and $\e$-entropy, as well as Temlyakov's
$(N,m)$-width.

On the positive side, we show that a very accurate {\it
non-linear} reconstruction is possible. It goes through a solution
of certain specific non-linear systems of algebraic equations. We
discuss the form of these systems and methods of their solution,
stressing their relation to Moment Theory and Complex Analysis.

Finally, we informally discuss two problems in Computer Imaging
which are parallel to the problems 1 and 2 above: compression of
still images and video-sequences on one side, and image
reconstruction from indirect measurement (for example, in Computer
Tomography), on the other.}}

\vspace{2 mm}
\begin{center}
------------------------------------------------
\vspace{2 mm}
\end{center}
This research was supported by the ISF, Grant No. 304/05, and by
the Minerva Foundation.

\end{titlepage}
\newpage

\leavevmode\kern0.2\hsize\begin{minipage}{0.8\hsize}\begin{small}Linear
problems are all linear alike; every non-linear problem is
non-linear in its own way.
\par
\medskip
\rightline{\emph{M. Livshitz}}\end{small}
\end{minipage}

\medskip

\section {Introduction}
\setcounter{equation}{0}

In this paper we discuss the following two basic problems:

1. How to represent functions with singularities (up to a
prescribed accuracy) in a compact way?

2. How to reconstruct such functions from a small number of
measurements?

We consider both the problems mainly from the point of view of a
comparison between linear and non-linear approaches.

We study in detail a model case of piecewise-constant functions on
$[0,1]$, which, as we believe, reflects many important issues of a
general situation. Considered as curves (or surfaces of higher
dimension) in the Hilbert space $L^2([0,1])$ the families of
piecewise-constant functions with variable jump-points form a nice
geometric object: so called ``crinkled arks". They are
characterized by the fact that any two their disjoint chords are
orthogonal. A remarkable classical fact is that any two such
curves are isometric, up to a scale factor. We reinterpret this
fact in a context of step-functions in one or two variables.

Next we study the problem of representability of
piecewise-constant functions by linear and semi-linear methods.
Our main tools in this problem are Kolmogorov's $n$-width and
$\e$-entropy (\cite{Kol1,Kol2}), as well as Temlyakov's
$(N,m)$-width (\cite{Tem2}). See also (\cite{Vet1,Vet3,Vet5}) and
references there for similar estimates.

Then we turn to the reconstruction problem. We start with a
negative result: based on our computation of Kolmogorov's
$n$-width of piecewise-constant functions, we provide limitations
on the accuracy of {\it linear} methods of reconstruction of such
functions from measurements.

On the contrary, we show, following \cite{Dav1, Dav2},
\cite{Mil1,Mil2,Mil3,Mil4,Mil5,Mil6},
\cite{Kis1,Kis2,Sar.Yom1,Sar.Yom2}, that a very accurate {\it
non-linear} reconstruction is possible. It goes through a solution
of certain specific non-linear systems of algebraic equations. We
discuss a typical form of these systems and certain approaches to
their solution, stressing the relations with Moment Theory and
Complex Analysis. See also \cite{Vet2,Vet4,Vet6} where a similar
approach is presented from a quite different point of view.

We believe that the key to a successful application of the
``algebraic reconstruction methods" presented in this paper to
real problems in Signal Processing lies in a ``model-based"
representation of signals and especially of images. This is a very
important and difficult problem by itself (see
\cite{Kun.Iko.Koc,Kou,Eld,Bri.Eli.Yom,Eli.Yom} and references
there). In the last section we informally discuss this problem
together with two other closely related problems in Computer
Imaging (which are parallel to the problems 1 and 2 above):
compression of still images and video-sequences on one side, and
image reconstruction from indirect measurement (for example, in
Computer Tomography), on the other.

\medskip

Our main conclusions are as follows:

\medskip

1. If we insist on approximating all the family of the
piecewise-constant functions, with variable positions of jumps,
{\it by the same linear subspace (Kolmogorov $n$-width)} then the
Fourier expansion is essentially optimal. Any other linear method
will provide roughly the same performance: with $n$ terms linear
combinations we get an approximation of order $1\over {\sqrt n}$.
This concerns both the ``compression" and the ``reconstruction
from measurements" problems.

\medskip

2. If for each individual piecewise-constant function we are
allowed to take {\it its own ``small" linear combination} of
elements of a certain fixed ``large" basis (``sparse
approximations") then with $n$ terms linear combination we get an
approximation of order $q^n, \ q<1$.

\medskip

3. The ``non-linear width" approach (Temlyakov's $(N,m)$-width)
provides a natural interpolation between the Fourier expansion,
the sparse approximations and the direct non-linear
representation.

\medskip

4. The ``naive" direct non-linear representation of
piecewise-constant functions, where we explicitly memorize the
positions of the jumps $0<x_i<1, \ i=1,\dots,N,$ and the values
$A_i$ of the function between the jumps, provides the best
possible compression (not a big surprise!). However, {\it these
parameters can be reconstructed from a small number of
measurements (Fourier coefficients) in a robust way, via solving
non-linear systems of algebraic equations}.

\medskip

5. Extended to piecewise-polynomials, and combined with a
polynomial approximation, the last result provides an approach to
an important and intensively studied problem of a noise-resistant
reconstruction of piecewise-smooth functions from their Fourier
data.

\medskip

Let us stress that the problem of an efficient reconstruction of
``simple" (``compressible") functions from a small number of
measurements has been recently addressed in a very convincing way
in the ``compressed sensing", ``compressive sampling", and
``greedy approximation" approaches (see
\cite{Can,Can.Don,Can.Rom.Tao,Don3,Don4,DeV,Tem1,Tem3} and
references there). Our approach is different, but some important
similarities can be found via the notion of ``semi-algebraic
complexity" (\cite{Yom1,Yom2}). We plan to present some results in
this direction separately.

\medskip

The third author would like to thank G. Henkin for very inspiring
discussions of some topics related to this paper. Both the
complexity of approximations and the moment inversion problem
intersect with Henkin's fields of interest, and we hope that some
of his results (see especially \cite{Hen1,Hen2,Hen3}) may turn out
to be directly relevant to the non-linear representation and
reconstruction problems discussed here.

\section{Families of piecewise-constant functions in $L^2([0,1])$}

In this paper we mostly concentrate on one specific case of a
piecewise-constant functions, namely, on the family of step (or
Heaviside) functions $H_t(x)$ defined on $[0,1]$ by $H_t(x)=1, \ x
\leq t$ and $H_t(x)=0, \ x > t$. All the results in Section 2
below remain valid (with minor modifications) for any family of
piecewise-constant functions on $[0,1]$ with a fixed number $N$ of
variable jumps.

A remarkable geometric fact about the curve $\mathcal{H}=\{H_t(x),
\ t \in [0,1]\} \subset L^2([0,1])$ is that any two its disjoint
chords are orthogonal. So the curve $\mathcal{H}$ changes
instantly its direction at each of its points: it is as
``non-straight" as possible. Such curves are called ``crinkled
arks" and we study them in more detail in Section 2.1.

Notice that a general family of piecewise-constant functions on
$[0,1]$ with a fixed number $N$ of variable jumps forms what can
be called a ``crinkled higher-dimensional surface" in
$L^2([0,1])$, at least with respect to the jump coordinates: any
two chords from the same point, corresponding to the jumps shifts
in opposite directions, are orthogonal.

\subsection{``Crinkled arcs"}

As above, we define the curve $\mathcal{H}:[0,1]\longrightarrow
L^2([0,1])$ by $\mathcal{H}_t=\chi_{[0,t]}.$ This curve is
continuous, and it has the following geometric property: any two
disjoint chords of it are orthogonal in $L^2([0,1])$. Indeed, such
chords are given by the characteristic functions of two
non-intersecting intervals. Intuitively, the curve $\mathcal{H}$
exhibits a ``very non-linear" behavior: its direction in
$L^2([0,1])$ rapidly changes.

Now let $\mathcal{X}$ be a general Hilbert space. \bd A curve
$\psi:[0,1]\longrightarrow {\mathcal{ X}}$ in a Hilbert space
${\mathcal{ X}}$ is called a crinkled arc if :
\begin{itemize}
\item it is continuous
\item any two disjoint chords of it are orthogonal, namely that
for $0\leq s<t\leq s'<t'\leq 1$ we have:
\end{itemize}
\be (\psi_t-\psi_s,\psi_{t'}-\psi_{s'})=0 \ee \ed More details are
given in the classical book of Halmos \cite{Hal}. See, in
particular, \cite[problems 5-6]{Hal}. The curve $\mathcal{H}$
provides the main example of a crinkled arc.

Crinkled curves are preserved by certain natural transformations.
Namely, one can perform
\begin{itemize}
\item  translation
\item scaling
\item reparametrization
\item application of a unitary operator.
\end{itemize}
Then the result would still be a crinkled arc. A simple and
surprising theorem is that these are the only possibilities to
obtain a crinkled arc, and any two arcs are connected by this
transformations: \bt Let $\psi:[0,1]\longrightarrow \mathcal{
X}_1$ and $\phi:[0,1]\longrightarrow {\mathcal{ X}}_2$ be two
crinkled arcs in two different Hilbert spaces. Then there are two
vectors $v_i\in{\mathcal{ X}}_i, \ i=1,2$, a reparametrization
$f:[0,1]\longrightarrow [0,1]$, a positive number $\alpha$ and a (partial)
isometry\footnote{a partial isometry $U:\mathcal{X}_1\rightarrow \mathcal{X}_2$
between Hilbert spaces is an isometry between $\text{Ker}U^\bot$ and
$\overline{\text{Im}U}$} $U:{\mathcal{ X}}_1\longrightarrow {\mathcal{ X}}_2.$ of
the Hilbert spaces s.t. \be U(\psi_{f(t)}-v_2)=\alpha\phi_t-v_1
\ee \et \pr See \cite[p.169]{Hal} \bc Let
$\psi:[0,1]\longrightarrow {\mathcal{ X}}$ be a crinkled arc. Then
it can be obtained from ${\mathcal{ H}}$ by a translation,
scaling, reparametrization, and an application of a unitary
operator between the appropriate Hilbert (sub-)spaces. \ec Therefore, if
we consider only geometric properties of curves inside the Hilbert
space then the curve $\mathcal{ H}$ can be taken as a model for
any crinkled ark.

While any two Hilbert spaces are isomorphic, their "functional"
realizations may be quite different. Consider, for example, the
space $L^2(Q^2)$ of the square integrable functions on the unit
two-dimensional cell $Q^2=[0,1]\times [0,1]$ (we shall later refer
to such functions as ``images").

\medskip

\begin{figure}[htbp]
   \begin{picture}(200,200)
   \put(0,150){$t_1$}
   \put(0,100){$t_2$}
   \put(0,50){$t_3$}

   \put(30,30){\framebox(80,40){\put(40,0){\circle*{15}}}}
   \put(30,80){\framebox(80,40){\put(40,0){\circle*{10}}}}
   \put(30,130){\framebox(80,40){\put(40,0){\circle*{5}}}}
   \put(130,30){\framebox(80,40){\thicklines\multiput(5,0)(1,0){35}{\line(0,1){20,5}}}}
   \put(130,80){\framebox(80,40){\thicklines\multiput(5,0)(1,0){25}{\line(0,1){20,5}}}}
   \put(130,130){\framebox(80,40){\thicklines\multiput(5,0)(1,0){15}{\line(0,1){20,5}}}}

   \end{picture}
   \caption{Two families of functions in $L^2(Q^2)$ that have the same Hilbert space
   geometric properties as $\mathcal{H} \subset L^2(Q^1)$ }
   \label{Pictures}
   \end{figure}
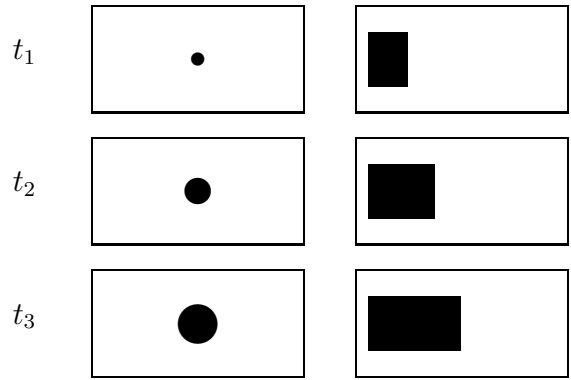

The two families of functions, which are shown in Figure
\ref{Pictures}, clearly represent crinkled arks in $L^2(Q^2)$.
Indeed, their disjoint chords are given by the characteristic
functions of certain concentric non-intersecting domains in $Q^2$,
and hence they are orthogonal. By Corollary 2.3, each of these
curves is isomorphic to the curve $\mathcal{H}$ in $L^2([0,1])$.
Let us state a general proposition in this direction. Consider a
family $D_t \subset Q^n, \ t\in [0,1],$ of ``expanding domains" in
the n-dimensional cell $Q^n=[0,1]^n$, $D_{t_1}\subset D_{t_2}$ for
any $t_1 < t_2$. Consider the curve $S(t)$ in $L^2(Q^n)$ defined
by $S(t)=\chi_{D_t}\in L^2(Q^n)$. \bp $S(t)$ is a crinkled
curve.\ep\pr Any two disjoint chords of the curve $S$ are given by
the characteristic functions of certain concentric
non-intersecting domains in $Q^n$, and hence they are orthogonal
in $L^2(Q^n)$.

\medskip

By Corollary 2.3, each of the curves $S$ obtained as above, is
isomorphic to the curve $\mathcal{H}$ in $L^2([0,1])$.

\begin{figure}[htbp]
\begin{pspicture}(13,6.5)

\psset{yunit=.01} \psset{xunit=0.1}
\pscustom[fillstyle=vlines*] {
\parabola(25,262.5)(70,600)
\parabola(115,262.5)(70,540)
}
\parabola(10,0)(70,600)
\parabola(10,46.7)(70,540)
\pscustom[fillstyle=hlines*] {
\parabola(20,241.7)(70,450)
\parabola(120,241.7)(70,405)
}
\parabola(10,150)(70,450)
\parabola(10,170)(70,405)
\end{pspicture}
\caption{Two chords of the family $C_t$}
   \label{Pictures2}
\end{figure}
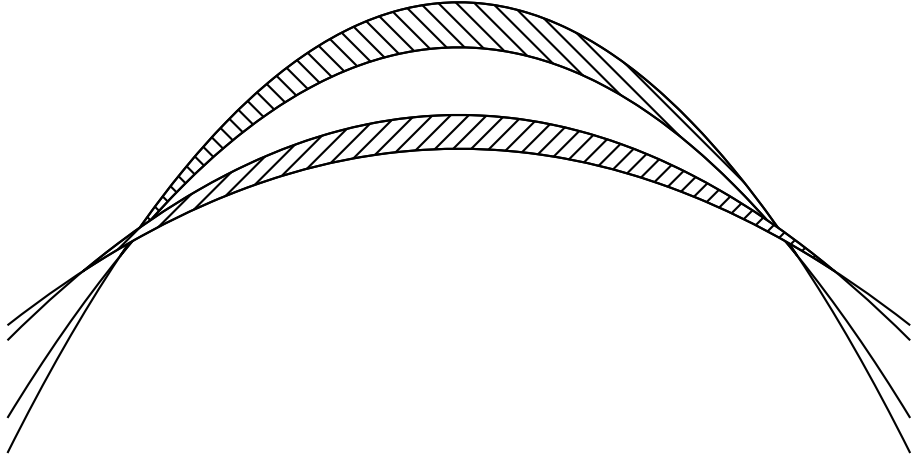

If the domains evolve in time in a more complicated way (in
particular, their boundaries are deformed in a non-rigid manner),
then the corresponding curve formed in $L^2(Q^n)$ may be not
exactly a crinkled arc. However, the following proposition shows
that {\it typically} such trajectories look like crinkled arcs
``in a small scale". \bp Let $C_t, \ t\in [0,1],$ be a generic
smooth family of closed non-intersecting curves in $Q^2$. Consider
a corresponding curve $C_t \subset {\mathcal {L}}$. Then the angle
between any two disjoint chords of $C_t$ tends to $\frac{\pi}{2}$
as these chords tend to the same point.\ep \pr Let us assume that
the curves $C_t(\tau)$ are parametrized by $\tau \in [0,1],\
C_t(0)=C_t(1)$. Because of the genericity assumption we can assume
that for each $t\in [0,1]$ the derivative ${\frac{\partial
C_t(\tau)} {\partial t}}$ has a finite number of zeroes
$\tau_1,\dots,\tau_m$ and it preserves its sign between these
zeroes. Therefore, the chords of $C_t$ are the characteristic
functions of the domains as shown on Figure \ref{Pictures2}.
Specifically, the intersections of these domains are concentrated
near the zeroes $\tau_1,\dots,\tau_m$ of ${\frac{\partial
C_t(\tau)} {\partial t}}$. Clearly, the area of the possible
overlapping parts of these domains is of a smaller order than the
area of the domains themselves.

\subsection{$\e$-Entropy of $\mathcal {H}$}
\setcounter{equation}{0}

From now on we compute the $\e$-entropy, the linear and non-linear
width only for the curve $\mathcal {H}$ in the space $L^2([0,1])$.
All these quantities depend only on the curve, and not on its
parametrization, and they are preserved by the isometries of the
ambient Hilbert spaces. To exclude the influence of the scalar
rescaling we can normalize our curves, for example, assuming that
the distance between the end-points is one. Then by Corollary 2.3
the $\e$-entropy, the linear and non-linear width are exactly the
same for each crinkled curve.

\medskip

Let us remind now a general definition of $\e$-entropy. Let $A
\subset X$ be a relatively compact subset in a metric space $X$.
\bd For $\e >0$ the covering number $M(\e,A)$ is the minimal
number of closed $\e$-balls in $X$ covering $A$. The binary
logarithm of the covering number, $H(\e,A)=\log M(\e,A)$ is called
the $\e$-entropy of $A$. \ed See \cite{Kol1,Kol.Tih} and many
other publications for computation of $\e$-entropy in many
important examples. Intuitively, $\e$-entropy of a set $A$ is the
minimal number of bits we need to memorize a specific element of
this set with the accuracy $\e$. Thus {\it it provides a lower
bound for the ``compression" of $A$, independently of the specific
compression method chosen}. \bp For the curve $\mathcal {H}$ in
the space $L^2([0,1])$ we have \be M(\e,\mathcal{H}) \asymp
(\frac{1} {\e})^2, \ H(\e,\mathcal{H}) \sim 2\log(\frac{1} {\e}).
\ee Here the sign $\asymp$ is used as an equivalent to the
inequality
$$ C_1 (\frac{1} {\e})^2 \leq M(\e,\mathcal{H}) \leq
C_2(\frac{1}{\e})^2$$ for certain $C_1$ and $C_2$, and for all
sufficiently small $\e$. The sign $\sim$ shows that $C_1$ and
$C_2$ tend to $1$ as $\e$ tends to zero. \ep \pr Let us subdivide
uniformly the interval $[0,1]$ into $N$ segments $\Delta_i$ by the
points $t_i={\frac i N}$. We have
$$\Vert \mathcal{H}(t_{i+1})-\mathcal{H}(t_{i})\Vert =
(\int_{t_i}^{t_{i+1}}dt)^{\frac1 2} = (\frac{1}{ N})^{\frac1 2}.$$
Hence for $\e={\frac1 2}({\frac1 N})^{\frac1 2}$ the $\e$-balls
covering different points $\mathcal{H}(t_i), \ i=1, \dots, N$ of
the curve $\mathcal{H}$ do not intersect. Thus, we need at least
$N$ such $\e$-balls to cover $\mathcal{H}$, while the $2\e$-balls
centered at the points $\mathcal{H}(t_i), \ i=1, \dots, N$ cover
the entire curve $\mathcal{ H}$. This completes the proof.

\subsection{Kolmogorov's $n$-width of $\mathcal{H}$}

Let $A \subset V$ be a centrally-symmetric set in a Banach space
$V$. \bd (\cite{Lor1,Kol2}). The Kolmogorov's $n$-width $W_n(A)$
of the set $A \subset V$ is defined as \be W_n(A) = \inf_{dim
L=n}\sup_{x\in A} \dist(x,L),
 \ee where the infinum is taken over all the
$n$-dimensional linear subspaces $L$ of $V$, and $\dist(x,L)$
denotes the distance of the point $x$ to $L$. \ed Intuitively,
$W_n(A)$ is the best possible approximation of $A$ by
$n$-dimensional linear subspaces of $V$. Let us define also
$N(\e,A)$ as the minimal $n$ for which $W_n(A)\leq \e$.

To make the Kolmogorov $n$-width comparable with the $\e$-entropy,
we define the notion of a {\it linear $\e$-entropy of $A$}, which
is the number of bits we need to memorize $A$ with the accuracy
$\e$, if we insist on a {\it linear approximation of $A$} (and if
we ``naively" memorize each of the coefficients in this linear
approximation): \bd A linear $\e$-entropy of $A$, $H_l(\e,A)$, is
defined by \be H_l(\e,A) = N(\e,A)\log (\frac{1}{ \e}).\ee \ed Now
we state the main result of this section: \bt \label{th:hest} For the curve
$\mathcal{H}$ in $L^2[0,1]$ we have $$\frac{1}{4\sqrt{n}} \leq
W_n(\mathcal{H})\leq \frac{2}{\pi\sqrt{n-1}}, \
N(\e,\mathcal{H})\asymp (\frac{1}{ \e})^2, \
H_l(\e,\mathcal{H})\asymp (\frac{1}{ \e})^2 \log(\frac{1}{ \e}).$$
\et \pr It is enough to prove the bounds for the $n$-width of
$\mathcal{H}$. The corresponding bound for $N(\e,\mathcal{H})$ and
$H_l(\e,\mathcal{H})$ follow immediately.

\medskip

Now, the upper bound for the $n$-width we obtain, considering the
Fourier series approximation of the Heaviside functions $H_t(x)$.
\be H_t(x)=\sum\limits_{k\in\mathbb{Z}}{a_ke^{2\pi i k x}} \ee
Then $a_0=t$ and $a_n=\frac{1-e^{-2\pi i n t}}{2\pi in}$ for
$n\neq 0$. We have $\vert a_n \vert \leq \frac{1}{\pi n}$. Hence
the $L^2$ error $f_n$ of the approximation of any $H_t$ by the
first $2n+1$ terms of its Fourier series satisfies

\be f_n \leq \bigl[\sum\limits_{m=n+1}^\infty{\frac{2} {\pi^2m^2}
}\bigr]^{\frac12} < \frac{\sqrt{2}}{\pi \sqrt{n}}. \ee And
therefore \be W_{2n+1}(\mathcal{H})\leq f_n \Longrightarrow
W_n(\mathcal{H})\leq \frac{2}{\pi \sqrt{n-1}}. \ee
\medskip

The proof of the lower bound we split into several steps. \bl For
a set $A_k=\{ e_i|(e_i,e_j)=\delta_{ij},1\leq i \leq k\}\subseteq L^2[0,1]$ and
$n<k$ the following inequality holds
\[W_n(A_k)\geq \sqrt{\frac{k-n}{k}}\]
\label{lemm:ess}. \el \pr Denote $W=\text{span}\{e_i\}$
and $P_W:L^2[0,1]\rightarrow L^2[0,1]$
the orthogonal projection on $W$. \\We take an n-dimensional
subspace $V$. We can assume that $V\subseteq W$ This is because
for $v\in V, \ a\in A_k$ we have:
\[||a-v||^2=||P_W(a-v)||^2+||(I-P_W)(a-v)||^2=\]
\[=||a-P_Wv||^2+||(I-P_W)v||^2\geq ||a-P_Wv||^2.\]
Therefore $\dist(A_k,V)\geq \dist(A_k,P_WV)$, and in order to
minimize the distance we can assume $V\subseteq W$. Denote
$P_V:W\longrightarrow W$ the orthogonal projection on $V$ in $W$. We
need to compute $\max_{1\leq i\leq k} \ ||(I-P_V)e_i||$. But
\be\max_{1\leq i\leq k} \ ||(I-P_V)e_i||\geq
\sqrt{\frac{1}{k}\sum\limits_{i=1}^k{||(I-P_V)e_i||^2}}.
\label{eas1} \ee On the other hand,
\be\sum\limits_{i=1}^k{||(I-P_V)e_i||^2}=
\sum\limits_{i=1}^k{((I-P_V)e_i,(I-P_V)e_i)}
=\sum\limits_{i=1}^k{((I-P_V)^2e_i,e_i)}=\label{eas2} \ee
\[=\sum\limits_{i=1}^k{((I-P_V)e_i,e_i)}=
\text{trace}_W({I-P_V})=k-n.\] The last equality is because
$I-P_V:W\rightarrow W$ is a projection into a $(k-n)$-subspace -
the orthogonal complement of $V$ in $W$. Combining equations
(\ref{eas1}),(\ref{eas2}) we have
\[\dist(A_k,V)\geq \sqrt{\frac{k-n}{k}}\]
\bc For any $d\in\mathbb{R}$,
\[W_n(dA_k)\geq |d|\sqrt{\frac{k-n}{k}}.\] \label{ort_wn_for}
\ec \bp\label{pr:hestlow} \[W_n(\mathcal{H})\geq \frac{1}{4\sqrt{n}}.\]
\ep \pr
Denote
$B_{k}=\{\chi_{\left(\frac{i-1}{k},\frac{i}{k}\right)}|1\leq i
\leq k\}$. This set is formed by $k$ orthogonal vectors of length
$\frac{1}{\sqrt{k}}$. Clearly
$\chi_{\left(\frac{i-1}{k},\frac{i}{k}\right)}= \mathcal{
H}_{\frac{i}{k}}-\mathcal{H}_{\frac{i-1}{k}}$, therefore
\[\dist(B_{k},V)\leq 2 \dist (\mathcal{H},V)\] for any vector space
$V$, and thus:
\[W_n(B_{k})\leq 2W_n(\mathcal{H}).\]
The norm of $\chi_{\left(\frac{i-1}{k},\frac{i}{k}\right)}$ is
$\frac{1} {\sqrt{k}}$, and according to Corollary \ref{ort_wn_for}
we have
\[W_n(\mathcal{H})\geq \frac{1}{2}W_n(B_{k})\geq \frac{1}{2}\frac{1}
{\sqrt{k}}\sqrt{\frac{k-n}{k}}.\] Taking $k=2n$ provides the
required result. This completes the proof of Proposition \ref{pr:hestlow} and
of Theorem \ref{th:hest}.

\subsection{Sparse representation of a step-function}

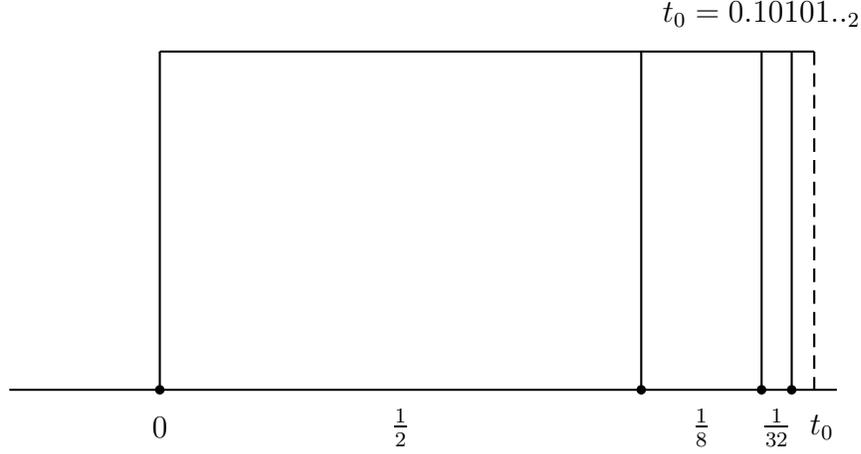
\begin{figure}[htbp]
\begin{pspicture}(11,7)
\rput(10,6){$t_0=0.10101.._2$} \psset{xunit=0.1}
\psline{-*}(0,1)(20,1) \psline{-*}(20,1)(84,1)
\psline{-*}(84,1)(100,1) \psline{-*}(100,1)(104,1)
\psline(104,1)(110,1)
\psline(20,5.5)(104,5.5) \psline(20,1)(20,5.5)
\psline(84,1)(84,5.5) \psline(100,1)(100,5.5)
\psline(104,1)(104,5.5)
\psline[linestyle=dashed](104,5.5)(107,5.5)
\psline[linestyle=dashed](107,5.5)(107,1) \rput(20,0.5){0}
\rput(52,0.5){$\frac{1}{2}$} \rput(92,0.5){$\frac{1}{8}$}
\rput(102,0.5){$\frac{1}{32}$} \rput(108,0.5){$t_0$}

\end{pspicture}
\caption{Approximation by Haar frame of $\chi_{t_0}$}
\label{pic:approx}
\end{figure}

Our main example of the family $\mathcal{H}$ of the step-functions
$H_t(x)$ allows us to illustrate also some important features of
``sparse representations".  Consider the Haar frame:
\[
\mathbf
{HF}=\Bigl\{\phi_{k,j}(x)=2^{k/2}\phi(2^k(x-j))|k\in\{0,1,2,\dotsc\},
j\in\{0,1,2,\dotsc, 2^k-1\}\Bigr\}
\]
where $\phi=\chi_{[0,1]}$. To get an approximation of a certain
fixed step-function $H_{t_0}(x)$ consider the binary
representation of $t_0$:
$$t_0=\sum^{\infty}_{r=1}{\frac{\alpha_r}{ 2^r}}, \
\alpha_r=0,1.$$ Then for each $n$ the sum \be \sum_{r\leq
n}\phi_{r,j_r}(x), \ j_r=\frac{1}{
2^r}\sum^{r-1}_{s=1}{\frac{\alpha_s}{ 2^s}}\ee leads to the
approximation of $H_{t_0}(x)$ in the Haar frame with the
$L^2$-error at most $(\frac1 2)^{\frac{n+1} 2}$. Indeed, the sum
in (6.1) is, in fact, a step-function $H_{t_1}$, with $t_1\leq
t_0$ and $t_0-t_1\leq (\frac1 2)^{n+1}$ (see Figure
\ref{pic:approx}).

\medskip

So to $\e$-approximate each {\it individual} step-function
$H_{t_0}$ via the Haar frame in the $L^2$-norm, we need only
$2\log(\frac{1} {\e})$ {\it nonzero terms in the linear
combination}. This provides a natural example of a ``sparse
representation".

\medskip

Notice, however, that if we fix the required approximation
accuracy $\e=(\frac1 2)^{\frac{n+1} 2}$, and then let the jump
point $t$ of $H_t$ change, then the elements of the Haar frame,
participating in the representation of different $H_t(x)$,
eventually cover all the $2^n$ binary step-functions of the $n$-th
scale. So altogether, to approximate the entire curve $\mathcal{H}
\subset L^2([0,1])$, we need the space of the dimension
$2^n=({\frac{1} {\e}})^2$. This agrees with the value of
$W_n(\mathcal{ H})$ computed above.

\subsection{$n$-term representation}

In order to quantify the ``sparsness" of different representations
(and, in particular, to include the previous example in a more
general framework) we call (following \cite[chapter 8]{DeV}) a
countable collection $\mathbf{D}$ of vectors in a Banach space a
dictionary, and define the error of the $n$-term approximation of
a single function $f$ by: \be \sigma_n(f,{\mathbf
D})=\inf_{w_i\in{\mathbf D},\alpha_i\in{\mathbb
C}}{\|f-\sum_{i=1}^n{\alpha_i w_i}\|}. \ee
We use three different dictionaries for $L_2[0,1]$: Fourier basis:
\[
\mathbf {FB}=\{e^{ikx}|k\in{\mathbb Z} \},
\]
Haar frame:
\[
\mathbf{
HF}=\Bigl\{\phi_{k,j}=2^{k/2}\phi(2^k(x-j))|k\in\{0,1,2,\dotsc\},
j\in\{0,1,2,\dotsc 2^k-1\}\Bigr\},
\]
and Haar basis:
\[\mathbf{
HB}=\Bigl\{\psi_{k,j}=2^{k/2}\psi(2^k(x-j))|k\in\{0,1,2,\dotsc\},
j\in\{0,1,2,\dotsc 2^k-1\}\Bigr\}\cup\{\phi\},
\]
where $\phi=\chi_{[0,1]}$ and $\psi=\chi_{[0,1/2]}-\chi_{[1/2,1]}.$\\
\medskip

Clearly,
\[
\sigma_n(H_t,\mathbf {FB})\asymp\frac{1}{n^{1/2}},
\]
which means, that the best $n$-term approximation in this case is
the same as the usual linear Fourier approximation. Also, we have
\[
\sigma_n(H_t,\mathbf {HF})\leq C 2^{-{\frac n 2}},
\]
\[
\sigma_n(H_t,\mathbf {HB})\leq C' 2^{-{\frac n 2}}.
\]
Remark: It is customary in the Approximation Theory to demand that
$n$-term approximation will be "computable"- so that it has
polynomial-depth search. This means that we can enumerate
our dictionary with a fixed enumeration $\mathbf
{D}=\{f_1,f_2,...\}$ in such a way that for a certain
polynomial $p:\mathbb {N}\rightarrow \mathbb {N}$ the $n$-terms of
approximation come from the first $p(n)$ terms of the dictionary,
see \cite{Can.Don}. Clearly, if we consider $\mathbf {HF}$ and
$\mathbf {HB}$ we will need to take each function from a different
level of Haar basis/frame, and therefore our search will have an
exponential depth.

\subsection{Temlyakov's non-linear width}

The following notion of a ``non-linear width" was introduced in
\cite{Tem2}: \bd Let $A$ be a symmetric subset in a Banach space
$\mathcal{ X}$. Then the $(N,m)$ width $W_{(N,m)}(A)$ is defined
as
\[W_{(N,m)}(A)=\inf_{{\mathcal{L}}_N\subseteq {\mathcal{L}}(X)_m, \
\vert{\mathcal{L}}_N\vert=N}{\sup_{f\in A}{\inf_{L\in\mathcal{
L}_N}{\dist(f,L)}}},\] where ${\mathcal{L}}(X)_m$ denotes the
collection of all the linear $m$-dimensional subspaces of
$\mathcal{ X}$. \ed The approximation procedure, suggested by this
notion, is as follows: given $N$ and $m$, we fix (in an optimal
way) a subset $\mathcal{ L}_N$ of $N$ different $m$-dimensional
linear subspaces $L_1,\dots,L_N$ in $\mathcal{ X}$. Then for each
specific function $f\in {\mathcal{X}}$ we first pick the most
suitable subspace $L_i$ in $\mathcal{ L}_N$, and then find the
best linear approximation of $f$ by the elements of $L_i$.

\medskip

The notion of a nonlinear width provides a ``bridge" between the
linear approximation and the approximation based on ``geometric
models". Indeed, ultimately the set $\mathcal{ L}_N$ may be just
the set formed by all the piecewise-constant functions (in our
main example), for all the values of the parameters, discrtetized
with the required accuracy. See Section 3 below where we analyze
in somewhat more detail this ``bridging" for the curve
$\mathcal{H}$.

\medskip
The set $\mathcal{ L}_N$ suggests a covering of the $A$ by $n$
sets
\[
V_i=\{g|\text{dist}(g,L_i)\leq \text{dist}(g,L_k);L_i,L_k\in \mathcal{ L}_N\}.
\]
Namely, the set $V_i$ contains the elements of $A$, that are best approximated
by the subspace $L_i$ from the collection $\mathcal{ L}_N$. In the next lemma,
we prove that we can replace $V_i$ by open sets.
\bl \label{WNm_char}
Let $\mathcal {O}_N$ denote the set of all the open covers
$\mathcal {U}=\{U_1,\dots,U_N\}$ of $A$ of cardinality $N$. Then
\[W_{(N,m)}(A)=\inf_{\mathcal {U}\in \mathcal {O}_N}
{\sup_{U_i\in \mathcal{U}}{W_m(U_i)}}.\] \el In other words, we
subdivide $A$ into $N$ open sets and check $m$-width on each of
the sets separately. Then the maximum $m$-width over $N$ sets is
the $(N,m)$-width of $A$.

\medskip

\pr Denote by $$d=W_{(N,m)}(A),$$ the left-hand side of the
equation and by $$e=\inf\limits_{\mathcal {U}\in \mathcal
{O},\#\mathcal {U}=N} {\sup\limits_{U_i\in
\mathcal{U}}{W_m(U_i)}},$$ it's right-hand side. For $\e>0$, we
interptret the definition of $d$ as existence of ${\mathcal{L}}_N$
a collection of $N$ $m$-dimensional subspaces of $X$ such that
\[\inf_{L_i\in{\mathcal{L}}_N,i=1..N}{\inf_{g\in L_i}\|f-g\|}
< d+\e\quad \forall f\in A\] Define $U_i=\{g|g\in
A,\|g-L_i\|<d+\e\}$. Clearly, $U_i$ are open in $A$ since the
distance is a continuous function. According to the definition of
$U_i$, $L_i$ approximates $U_i$ with accuracy $d+\e$ and therefore
$W_m(U_i)\leq d+\e$. We conclude that \be e\leq \sup_{i}
W_m(U_i)\leq d+\e. \label{eld} \ee In the other direction, let
$\bigcup\limits_{i=1}^{N} U_i=A$ such that $W_m(U_i)<e+\e$. For each
$i$ we find $L_i$ an $m$-dimensional subspace s.t. $$
\dist(g,L_i)<e+2\e,\quad\forall g\in U_i, $$ Then form
${\mathcal{L}}_N=\{L_1,..,L_N\}.$ Clearly,
\[\sup_{f\in A}{\inf_{L\in{\mathcal{L}}_N}{\dist(f,L)}}<e+2\e\]
and therefore \be d<e+2\e \label{dle} .\ee Taking
$\e\longrightarrow 0$ in the inequalities (\ref{eld}),
(\ref{dle}), we get the required equality.
\medskip
\\ In what follows, we take $\mathcal{X}=L_2([0,1])$.
\bp
\[W_{(N,m)}(\mathcal{H}) \asymp \frac{1}{\sqrt{Nm}}.\]
\label {pr:nwidth} \ep \pr
\\
To establish an upper bound, define
$$L_k=\text{span}\{\chi_{[\frac{k-1}{N}+\frac{n}{Nm},
\frac{k-1}{N}+\frac{n+1}{Nm}]}:0\leq n\leq m-1\},\quad k=1..N$$ Each $L_k$
approximates $\mathcal{H}\bigl|_{[\frac {k-1}{N},\frac{k}{N}]}$
within an error of $\frac{1}{\sqrt{2Nm}}$. Therefore $W_{(N,m)}$
is bounded above by an error of $\{L_k \}_{k=1}^{N}$. \\
In order to establish the lower bound, we prove a variant of a
Lemma \ref{lemm:ess}. \bl For a set $A_k=\{
e_i|(e_i,e_j)=\lambda_i^2\delta_{ij},1\leq i \leq k,\lambda_i>0\}$
and $n<k$ the following inequality holds
\[W_n(A_k)\geq \sqrt{\frac{k-n}{k}}\min\limits_i{\lambda_i}.\]
\label{lemm:mod} \el The difference with Lemma \ref{lemm:ess} is
that we allow
orthogonal vectors with varying lengths.\\
\pr Let $V$ be a $n$-dimensional space and $W=\text{span}\{e_i\}$.
Just like in Lemma \ref{lemm:ess}, we can assume $V\subseteq W$ .
Denoting the orthogonal projection of $V$ on $W$ by $P$, we are
required to compute $\max\limits_{i}\|(I-P)e_i \|$. \be
\max\limits_{i}\|(I-P)e_i \|\geq
\sqrt{\frac1k\sum\limits_{i=1}^k\|(I-P)e_i \|^2}\geq
\min\limits_{i} {\lambda_i}
\sqrt{\frac1k\sum\limits_{i=1}^k{\|(I-P)\frac{e_i}{\lambda_i}
\|^2}}. \label{eq:lm1} \ee Since $\frac{e_i}{\lambda_i}$ are
orthonormal then \be
\sum\limits_{i=1}^k{\|(I-P)\frac{e_i}{\lambda_i}\|^2}=
\sum\limits_{i=1}^k{((I-P)\frac{e_i}{\lambda_i},\frac{e_i}
{\lambda_i})}=\text{trace}(I-P)=k-n. \label{eq:lm2} \ee Combining
equations \ref{eq:lm1} and \ref{eq:lm2} we get the required
result.
\bigskip
\\We return to the proof of the Proposition \ref{pr:nwidth}. We
employ Lemma \ref{WNm_char}. Let $\{U_i\}_{i=1..n}$ be an open
cover of $\mathcal{H}$. Define
\[V_i=\mathcal{H}^{-1}(U_i).\]
Namely, $V_i\subseteq [0,1]$ contains all the $t$'s such that
$\mathcal{H}(t)\in U_i$. $V_i$ are open in $[0,1]$, since
$\mathcal{H}$ is continuous. The collection $\{V_i\}_{i=1..n}$ is
a covering of $[0,1]$ since $\{U_i\}$ is a cover of $\mathcal{H}$.
Because $V_i$ is open, we can find $2m+1$ points $\lambda_k$ in
$V_i$ such that $\lambda_k-\lambda_{k+1}\geq
\frac{\text{meas}(V_i)}{2m}-\e$, for any $\e>0$, where
$\text{meas}(V_i)$ denotes here the Lebesgue measure of $V_i$. We
apply Lemma \ref{lemm:mod} to
$B_{2m}=\{\chi_{[\lambda_k,\lambda_{k+1}]}:1\leq k \leq 2m \}$. Since
$\|\chi_{[\lambda_k,\lambda_{k+1}]}\|=\sqrt{\lambda_{k+1}-\lambda_k}\geq
\sqrt{\frac{\text{meas}(V_i)}{4m}}-\e$, the application of Lemma
\ref{lemm:mod} gives
\[
W_m(B_{2m})\geq  \sqrt{\frac{\text{meas}(V_i)}{4m}}-\e.
\]
But $\chi_{[\lambda_k,\lambda_{k+1}]}=\mathcal{H}_{\lambda_{k+1}}-
\mathcal{H}_{\lambda_{k}}$.Denote
$S=\{\mathcal{H}_{\lambda_k}:1\leq k \leq 2m+1\}\subseteq U_i$. For
any vector space $V$, we have
\[ \dist(B_{2m},V)\leq 2 \dist (S,V),\]
and therefore \be W_m(U_i)\geq W_m(S)\geq  \frac 12 W_m(B_{2m})
\geq \frac14\sqrt{\frac{\text{meas}(V_i)}{m}}-\e. \ee But since
$V_i$ cover $[0,1]$ , we have
$\sum\limits_{i=1}^N{\text{meas}(V_i)}\geq 1$ and so
$$\max\limits_i\text{meas}(V_i)\geq \frac{1}{N}.$$
Therefore \be \max_i W_m(U_i)\geq \frac{1} {4\sqrt{Nm}}-\e, \ee
for any open cover of $\mathcal{H}$. And so, according to Lemma
\ref{WNm_char} \be W_{(N,m)}(\mathcal{H})\geq \frac{1}
{4\sqrt{Nm}}-\e \label{eq:wnm_int}. \ee Thus we obtain the
required lower bound after we take
 $\e\longrightarrow
0$.

\section{Linear versus Non-Linear Compression: some conclusions}

In this section we summarize the above results, interpreting them
as the estimates of the``compression" of the family $\mathcal{H}$
(and of other families of piecewise-constant functions): {\it how
many bits do we need to memorize an arbitrary jump-function $H_t$
in $\mathcal{H}$ with the $L^2$-error at most $\e$, via different
representation methods?}

\subsection{$\e$-entropy}

Let us start with the $\e$-entropy: by Proposition 2.7,
$H(\e,\mathcal{H})\sim 2\log(\frac{1}{ \e})$. This is the lower
bound on the number of bits in any compression method.

\subsection{``Model-based compression"}

Let us consider a ``non-linear model-based compression" which in
the case of the jump-functions takes an extremely simple form: we
use the ``library model" $H_t(x)$ to represent itself, and we
memorize just the specific value of the parameter $t$. Quite
expectedly, ``compression" with this model requires exactly the
number of bits prescribed by the $\e$-entropy. Indeed, since the
$L^2$-norm of $H_{t_2}-H_{t_1}$ is $\sqrt{t_2-t_1},$ we have to
memorize $t$ with the accuracy $\e^2$. This requires exactly
$2\log(\frac{1}{\e})$ bits.

\subsection{``Linear" compression}

Let us assume now that, given the required accuracy $\e$, we
insist on a representation of the functions $H_t(x)$ in a fixed
basis, the same for each $t$. On the other hand, we allow the
approximating linear space to depend on $\e$. This leads to the
Kolmogorov $n$-width, as defined in Section 2.3. We store each
coefficient with the maximal error $\e$, so we allow for it
$\log(\frac{1}{ \e})$ bits (and thus we ignore a very special
``sparse" nature of the representation of $H_t(x)$ in some special
bases, for instance, in the Haar frame, discussed in Section 2.4).
Then the number of bits required is given by the ``linear
$\e$-entropy" $H_l(\e,\mathcal{H})$, introduced in Section 2.3. By
Theorem 2.10, we have
$$H_l(\e,\mathcal{H})\asymp (\frac{1}{ \e})^2 \log(\frac{1}{ \e}).$$ In
fact, to get a representation with this amount of information
stored, we do not need all the freedom provided by the definition
of $n$-width. It is enough to fix the approximating space to be
the space of trigonometric polynomial for any required accuracy
$\e$. Then to approximate $H_t$ with the $L^2$-accuracy $\e$ we
take the Fourier polynomial $F^n_t$ of $H_t$ of degree
$n=\frac{1}{ \e^2}$ and memorize its coefficients with the
accuracy $\frac{\e}{n}$.

\subsection{``Non-linear width" compression}

In \cite{Tem2} a notion of a ``non-linear $(N,m)$-width" has been
introduced (see Section 2.6 above). It suggests the following
procedure for approximating functions $H_t(x)$: given the required
accuracy $\e$, we fix a subset of $N$ $m$-dimensional linear
subspaces $L_1,\dots,L_N$ in $L^2[0,1]$. Then for each specific
function $H_t(x)$ we first pick one of the subspaces $L_i$ (the
most suitable), then find the best linear approximation of
$H_t(x)$ by the elements of $L_i$, and finally memorize the
coefficients of the best linear approximation found.

Let us estimate the number of bits required in this approach. By
Proposition 2.16, for the non-linear $N,m$-width of $\mathcal{H}$
we have
$$W_{(N,m)}(\mathcal{H}) \asymp \frac{1}{\sqrt{Nm}}.$$ Given the required
accuracy $\e$, we have to fix the parameters $N$ and $m$ in such a
way that $\frac{1}{\sqrt{Nm}}\leq \e.$ Therefore, for each choice
of $m$ between $1$ and $\left(\frac{1}{ \e}\right)^2$ we have to
take $N=\frac{1}{m}\left(\frac{1}{\e}\right)^2$. To memorize the
choice of the space $L_i$ we need then $\log
N=2\log(\frac{1}{\e})-\log \ m$ bits. To memorize the coefficients
we need $m \log(\frac{1}{ \e})$ bits. Hence, the total amount of
bits is
$$(m+2)\log(\frac{1}{\e})-\log \ m. $$ Certainly, the best choice is
$m=1$: we just take $N=(\frac{1}{\e})^2$ elements $H_{t_i}, \
t_i=\frac{i}{N},$ and approximate $H_t$ with the nearest among
$H_{t_i}$. This is, essentially, the same as the ``model-based"
representation in Section 3.2 above.

\subsection{``Sparse" representation}

Till now the comparison was in favor of a model-based approach.
Let us consider now the Haar frame representation of $H_t(x)$
considered in Section 2.4 above. This is the most natural
competitor, both because of its theoretical efficiency, and since
many modern practical approximation schemes are based on sparsness
considerations (see \cite{Can,Can.Don,Tem1,Tem2,Tem3}).

By the computation of Section 2.4, to approximate each {\it
individual} step-function $H_{t_0}$ via the Haar frame in the
$L^2$-norm, we need only $m=2\log(\frac{1}{\e})$ of the {\it
nonzero terms in the linear combination}. Moreover, each
coefficient in this linear combination is $1$. So to memorize
$H_{t_0}$ via the Haar frame it is enough to specify the position
of $m=2\log(\frac{1}{\e})$ nonzero elements among the total Haar
frame of cardinality $2^m=({\frac{1}{\e}})^2$. We need
$$\log{\frac{(2^m)!} {m!(2^m-m)!}}\asymp m^2\asymp
[\log(\frac{1}{\e})]^2$$ bits to do this.

We get a little bit more information to store than in the
``model-based" approach. Also, it may look not natural to
approximate such a simple pattern as a jump of a step-function
with a geometric sum of shrinking signals. However, the main
problem is that if we let the jump point $t$ of $H_t$ change, then
the elements of the Haar frame, participating in the
representation of different $H_t(x)$, jump themselves in a very
sporadic way, and eventually cover all the $2^m$ binary Haar frame
functions of the $m$-th scale.

Notice also that from the point of view of the non-linear width
(Section 2.6 above) the considered Haar frame representation takes
an intermediate position: here $m=2\log(\frac{1}{\e})$. But any
subspace
$L=\text{span}\{\chi_{[t_{k_i},t_{k_i}+2^{-k_i}]:i=1..m}\}$ can
cover only $\mathcal{H}_{t_{k_i}}$ and their $\e$-neighborhoods
and therefore $L$ covers with the accuracy $\e$ only a set of
measure $4\e^2(\log(\frac{1}{\e})+1)$ out of the entire interval
$[0,1]$ of parameters. Thus to cover the entire interval we will
need $N\asymp \frac{1}{\e^2\log(\frac{1}{\e})}$ subspaces. We
conclude that  $\e\asymp\frac{1}{\sqrt{Nm}}$, in agreement with
Proposition \ref{pr:nwidth}. The required number of bits is
\[
\left(\log\left(\frac{1}{\e}\right)\right)^2+2\log\left(\frac{1}
{\e}\right)-\log\log\left(\frac{1}{\e}\right)+\log 2.
\]

\medskip

So it would be much more natural and efficient to represent a
``video-sequence" $\mathcal{H}=\{H_t(x), \ t\in [0,1]\}$ by a
moving model than to follow the jumping parameters in a sparse
Haar representation for variable $t$. This conclusion certainly is
not original. The problem is to get a full quality model-based
geometric representation of real life images and video-sequences!

\section{Non-linear Fourier inversion}

Now we turn to our second main problem: how to reconstruct
functions with singularities (piecewise-constant functions) from a
small number of measurements? Let us assume that our
``measurements" are just the scalar products of the function $f$
to be reconstructed with a certain sequence of basis functions. In
particular, below we assume our measurements to be the the Fourier
coefficients of $f$ or its moments. This is a realistic assumption
in many practical problems, like Computer Tomography.

The rate of Fourier approximation of a given function and the
accuracy of its reconstruction from partial Fourier data is
determined by regularity of this function. For functions with
singularities, even very simple, like the Heaviside function, the
convergence of the Fourier series is very slow. Hence a
straightforward reconstruction of the original function from its
partial Fourier data (i.e. forming partial sums of the Fourier
series) in this cases is difficult. It also involves some
systematic errors (like the so-called Gibbs effect).

\medskip

Let us show that no {\it linear} reconstruction method can do
significantly better that the straightforward Fourier expansion.

\bt Let the function acquisition process comprise taking $n$
measurements (linear or non-linear) $m_i(f), \ i=1,\dots,n$  of
the function $f$, together with a consequents processing $P$ of
these measurements. If the processing operator $\hat f =
P(m_1,\dots,m_n)$ is a linear operator from ${\mathbb R}^n$ to
$L^2([0,1])$ then for some $f \in \mathcal H$ the error $\vert
\vert f - \hat f \vert \vert$ is at least $C_1{1\over {\sqrt
n}}$.\et \pr This follows directly from Theorem 2.10 above.
Indeed, the $n$-dimensional linear subspace $Im(P)$ cannot
approximate all the functions in $\mathcal H$ with the error
better than the Kolmogorov $n$-width of $\mathcal H$.

\medskip

If we have no a priori information on $f \in L^2([0,1])$ then
probably the straightforward Fourier reconstruction as above
remains the best solution. However, in our case we know that $f$
is a piecewise constant function. It is completely defined by the
positions of its jumps $0<x_i<1, \ i=1,\dots,N$ and by its values
$A_i$ between the jumps. So let us consider $x_i$ and $A_i$ as
unknowns and let us substitute these unknowns into the integral
expression for the Fourier coefficients. We get certain analytic
expressions in $x_i$ and $A_i$. Equating these expressions to the
measured values of the corresponding Fourier coefficients we get a
system of nonlinear equations on the unknowns $x_i$ and $A_i$. Let
us write down this system explicitly.

\subsection{Fourier inversion system}

Let $f(x)=\sum_{-\infty}^{\infty} c_k e^{2\pi i k x}$ be the
Fourier expansion of $f$. Here $c_k={1\over {2\pi i}}\int_0^1
f(t)e^{-2\pi i k t}dt.$ Taking into account a special form of $f$
as given above we obtain $c_k={1\over {2\pi i k}}
[-A_0+\sum_{i=1}^N(A_{i-1}-A_i)e^{-2\pi i k x_i}+A_N e^{-2\pi i
k}]$. Here $A_0$ is the value of $f$ on the leftmost continuity
interval. Denoting $-{2\pi i k}c_k$ by $\hat c_k$ and $e^{-2\pi i
x_i}$ by $z_i$, we finally get the following infinite system \be
A_0+\sum_{i=1}^N(A_i-A_{i-1})z_i^k - A_N e^{-2\pi i k} = \hat c_k,
\ k \in \mathbb Z.\ee The unknowns in system (4.1) are $A_j, \
j=0,\dots,N$ which enter this system in a linear way, and $z_i, \
i=1,\dots,N,$ entering it non-linearly.

System (4.1) classically appears in Pade Approximation. Very
similar systems appear in a reconstruction of plane polygonal
domains from their moments (\cite{Mil1,Mil2,Mil3,Mil4,Mil5,Mil6}).
A detailed investigation of a larger class of systems similar to
(4.1) is given in \cite{Kis1,Kis2}. In particular,  we have the
following result: \bt Assume that $c_k$ in the right-hand side of
(4.1) are Fourier coefficients of a piecewise-constant function
$f$ with $A_i\neq A_{i+1}, \ i=0,\dots,N.$ Then each subsystem of
(2.1) obtained by taking from it certain $2N+1$ subsequent
equations has a unique solution $\{A_j, \ j=0,\dots,N\}, \ \{z_i,
\ i=1,\dots,N\}$, with $x_i={1\over {-2\pi i}}\log {z_i}$ being
the jump points of $f$ and $A_j$ being the values of $f$ on its
continuity interval.\et We give a sketch of the proof, following
\cite{Kis1,Kis2}, in Section 4.2 below.

\smallskip

Thus solving an appropriate subsystem of system (4.1) we find the
jumps and the intermediate values of $f$, so we reconstruct $f$
exactly. If $f$ had $N$ jumps we need only $2N+1$ Fourier
coefficients to reconstruct it.

\subsection{Other examples of the inversion systems}

Let us start with another system which essentially coincides with
(4.1). To simplify the presentation we shall consider instead of
the Fourier coefficients of the function $g(x), \ x \in [0,1]$ the
moments $m_k(g)=\int_0^1 x^kg(x)dx.$

\subsubsection{Linear combination of $\delta$-functions}

Let $g(x)=\Sigma_{i=1}^n A_i\delta(x-x_i)$. For this function we
have \be m_k(g)=\int_0^1 x^k \Sigma_{i=1}^n A_i\delta(x-x_i) dx =
\Sigma_{i=1}^n A_i x_i^k.\ee So assuming that we know the moments
$m_k(g)=\alpha_k, \ k=1,\dots, 2n-1,$ we obtain the following
non-linear system of equations for the parameters $A_i$ and $x_i,
\ i=1,\dots,n,$ of the function $g$:

\begin{eqnarray}
\Sigma_{i=1}^n A_i = \alpha_0, \nonumber \\ \Sigma_{i=1}^n A_i x_i
= \alpha_1, \nonumber \\ \Sigma_{i=1}^n A_i x_i^2 = \alpha_2,
\nonumber \\ .................... \nonumber \\ \Sigma_{i=1}^n A_i
x_i^{2n-1} = \alpha_{2n-1}.
\end{eqnarray} This system can be
solved as follows: consider the moments generating function
$$I(z)=\sum_{k=0}^{\infty} m_k(g) z^k$$. The representation (4.2)
of the moments immediately implies that \be I(z)= \Sigma_{i=1}^n
{A_i\over {1-zx_i}}.\ee So it remains to find explicitly the
rational function $I(z)$ from the first $2n$ its Taylor
coefficients $\alpha_0,\dots,\alpha_{2n-1}$.

\medskip

To do this we remind that the Taylor coefficients of a rational
function satisfy a linear recurrence relation of the form \be
m_{r+n}=\Sigma_{j=0}^{n-1} C_j m_{r+j}, \ r=0,1,\dots.\ee Since we
know the first $2n$ Taylor coefficients
$\alpha_0,\dots,\alpha_{2n-1}$, we can write a {\it linear} system
on the unknown recursion coefficients $C_l$:
\begin{eqnarray} \Sigma_{j=0}^{n-1} C_j \alpha_{j}=\alpha_{n},
\nonumber \\
\Sigma_{j=0}^{n-1} C_j \alpha_{j+1}=\alpha_{n+1}, \nonumber \\
\Sigma_{j=0}^{n-1} C_j \alpha_{j+2}=\alpha_{n+2}, \nonumber \\
................., \nonumber \\
\Sigma_{j=0}^n C_j \alpha_{j+n}=\alpha_{2n-1}.\end{eqnarray}
Solving linear system (4.6) with respect to the recurrence
coefficients $C_j$ we find them explicitly. For a solvability of
(4.6) see \cite{Nik.Sor,Kis1,Kis2,Sar.Yom1,Sar.Yom2}. See also
\cite{Mil2,Mil3,Mil4,Mil5}.

Now the recurrence relation (4.5) with known coefficients $C_l$
and known initial moments allows us to easily reconstruct the
generating function $I(z)$ and hence to solve (4.3).

\subsubsection{Algebraic functions}

Let now $g(x)$ be an algebraic function on $[0,1]$. By definition,
$y=g(x)$ satisfies an equation \be
a_n(x)y^n+a_{n-1}(x)y^{n-1}+\dots + a_1(x)y+a_0(x)=0,\ee where
$a_n(x),\dots,a_0(x)$ are polynomials in $x$ of degree $m$.
$d=m+n$ is, by definition, the degree of $g$.

A general method for the non-linear inversion of the moment
(Fourier) transforms of algebraic functions is given in
\cite{Kis2}. Its ``quantitative form is given in \cite{Sar.Yom2}.
Here we analyze only one special case. Assume that the algebraic
curve $y=g(x)$ is a rational one. This means that it allows for a
rational parametrization \be x=P(t), \ y=Q(t).\ee The moments
$m_k(g)$ given by $ m_k(g)=\int_0^1 x^k g(x)dx, \ k=0,1,\dots,$
now can be expressed as \be m_k(g)=\int_0^1 P^k(t)Q(t)p(t)dt,\ee
where $p$ denotes the derivative $P'$ of $P$. Moments of this form
naturally appear in a relation with some classical problems in
Qualitative Theory of ODE's - see \cite{bfy1,bfy2},
\cite{by,chr,yp}.

\medskip

Our problem can be reformulated now as the problem of explicitly
finding $P$ and $Q$ from knowing a certain number of the moments
$m_k$ in (4.9). Of course, in general we cannot expect this system
of nonlinear equations to have a unique solution. Indeed, while
the function $y=g(x)$ is determined by its moments in a unique
way, the {\it rational parametrization} of this curve in general
is not unique. In particular, let $W(t)$ be a rational function
satisfying $W(0)=0, \ W(1)=1$. Substituting $W(t)$ into $P$ and
$Q$ we get another rational parametrization of our curve: $x=\hat
P(t), \ y=\hat Q(t)$ with $\hat P(t)=P(W(t)), \ \hat
Q(t)=Q(W(t))$. Consequently, the ``inversion problem" for system
(4.9) is:

\medskip

{\it To characterize all the solutions of system (4.9) and to
provide an effective way to find these solutions}.

\medskip

A special case of the inversion problem is the ``Moment vanishing
problem":

\medskip

{\it To characterize all the pairs $P,Q$ for which the moments
$m_k$ defined by (4.9) vanish}.

\medskip

In spite of a very classical setting (we ask for conditions of
orthogonality of $Q$ to all the powers of $P$!) this problem has
been solved only very recently (\cite{pak2}). It plays a central
role in study of the center conditions for the Abel differential
equation (see \cite{bfy1,bfy2}, \cite{by,chr,yp},
\cite{pak1,pak2}.

\subsubsection{Functions of two variables}

The approach to reconstruction of piecewise-smooth
(piecewise-polynomial) functions of one variable discussed above
can be extended to two (and more) variables. The case of
characteristic functions of polygonal plane domains is considered
in \cite{Mil2,Mil3,Mil4,Mil5}. Some initial instances of the
reconstruction problem of piecewise-polynomial functions of two
variables are considered in \cite{Sar.Yom1}. Even the most initial
examples in two dimensions provide an exciting variety of
non-linear system bringing us to the very heart of Analysis. Let
us mention here only one example and a few of the most directly
related references.

We want to reconstruct a function $f(x,y)$ of two variables from
its moments \be m_{kl}(f)=\int\int x^k y^l f(x,y) dx dy. \
k,l=0,1,\dots.\ee Assume that $f$ is a $\delta$-function along a
rational curve $S$, i.e. for any $\psi(x,y)$ we have $\int\int
f\psi dx dy = \int_{S}\psi(x,y)dx$.

Let \be x=P(t), \ y=Q(t), \ t \in [0,1] \ee be a rational
parametrization of $S$. The moments now can be expressed as \be
m_{kl}(f)=\int_0^1 P^k(t)Q^l(t)p(t)dt,\ee where $p$ denotes the
derivative $P'$ of $P$. The study of the double moments of this
form bring us naturally to the recent work of G. Henkin
\cite{Hen1,Hen2,Hen3}. Indeed, the vanishing condition for the
moments (4.12) is given by Wermer's theorem (\cite{Wer}):
$m_{kl}(f)\equiv 0$ if and only if $S$ bounds a complex $2$-chain
in ${\mathcal C}^2$. In general, if the moments $m_{kl}(f)$ do not
vanish identically, then the local germ of complex analytic curve
$\hat S$ generated by $S$ in ${\mathcal C}^2$ does not ``close up"
inside ${\mathcal C}^2$. G. Henkin's work (\cite{Hen1,Hen2,Hen3}),
in particular, analyzes various possibilities of this sort in
terms of the ``moments generating function". We expect that a
proper interpretation of the results of \cite{Hen1,Hen2,Hen3} can
help also in understanding of the moment inversion problem.

\medskip

There are many other results closely related to our problem (see
references in \cite{Mil2,Mil3,Mil4,Mil5}, \cite{Sar.Yom1,
Sar.Yom2}. Here we mention in addition only \cite{Put1,Put2}
where, in particular, the problem of a reconstruction of plane
``quadrature domains" from their double moments is considered, and
results on moments on Semi-Algebraic sets and positivity (see
\cite{Kuh1,Kuh2,Sch} and references there).

\subsection{Robustness of solutions of (2.1)}

Let us return to functions of one variable. The assumption of $f$
being a piecewise-constant function may look too unrealistic in
applications. However, the methods can be extended to
piecewise-polynomial and ultimately to piecewise-smooth functions
(see \cite{Kis1,Kis2,Sar.Yom1,Sar.Yom2}). The last class is of
major importance in applied Analysis and Signal Processing, and
the problem of a reconstruction of such functions from their
measurements (Fourier data) is at present actively investigated
(see \cite{Tad1}-\cite{Tad4}, \cite{Vet1,Vet3,Vet5} and references
there).

The key issue in the extension of the above methods to
piecewise-smooth functions is a robustness of solutions of (4.1),
(4.3), (4.9) and similar systems. In particular, what happens if
we take more than exactly $2N+1$ consequent equations in (4.1),
and because of the noise in our measurements the right hand side
is not exactly a sequence of the Fourier coefficients of a
piecewise-constant function? Some important results in this
direction can be found in \cite{Mil2,Mil3,Mil4,Mil5}.

We further investigate these problems in \cite{Sar.Yom2}. Our
initial considerations show that one can define a robust procedure
for solving systems like (4.1),(4.3) for any right-hand side,
taking more equations than $2N+1$ and replacing the exact solution
by the least-square fitting. Notice, however, that we apply this
procedure not to the original non-linear system (say, (4.3)) but
to a linear system (4.6) for the parameters $C_j$ of the linear
recurrence relation, satisfied by the Fourier coefficients
(moments) $c_k$ ($m_k$) of any piecewise-constant function.

We expect also that taking more than the minimal number of the
measurements, and hance of the equations in (4.1) (say, twice the
minimal number) can strongly improve the robustness of the
solution. This conclusion is supported by recent results in
\cite{Zah,Yom.Zah} where we investigate similar problems for
Hermite interpolation and Hermite least-square fitting.

\subsection{Piecewise-smooth functions}

We expect that applying the above results to the case of
piecewise-smooth functions we can get, in particular, the
following result:

\medskip

\noindent{\bf Conjecture.} {\it Let $f$ be a piecewise $C^k$
function on $[0,1]$ with $N$ discontinuity points $x_i$. Assume
that the $C^k$-norm of $f$ on each continuity interval does not
exceed $M$. Assume also that a distance $\vert x_i-x_j \vert \geq
D$ for $i \neq j$, and that a jump of $f$ at each of its
discontinuity points $x_i$ is at least $J$. Then for each $n>2N+1$
the points $x_i$ and the values of $f$ between the points $x_i$
can be reconstructed from the first $n$ Fourier coefficients of
$f$ with the accuracy $C\over {N^k}$, where the constant $C$
depends only on $M, \ N, \ J$ and $D$.}

\medskip

We expect also that $x_i$ and the approximating polynomials of $f$
on its continuity intervals are provided by {\it universal
analytic expressions in $c_k$} (see \cite{Kis1,Kis2}). In
\cite{Sar.Yom2} we prove a weaker version of this result (with a
weaker estimate of the approximation accuracy). The main steps in
the proof are as follows:

\medskip

1. We fix an approximation accuracy $\e>0$. We approximate $f$ up
to $\e$ by a piecewise-polynomial $\Delta P$ of degree $d$ on each
of its continuity intervals. By classical Approximation Theory
this can be achieved with $d=C_1({1\over {\e}})^{1\over k}$.

\medskip

2. We consider the jump points $x_i$ and all the coefficients of
the piecewise polynomials constituting $\Delta P$ as the unknowns,
and substitute them into a system (*) similar to (4.1), which is
constructed ``once forever" for piecewise-polynomials of degree
$d$ (see \cite{Kis1,Kis2,Sar.Yom1,Sar.Yom2}). As the right-hand
side we take the Fourier coefficients $c_k$ of $f$. By the choice
of $\Delta P$ we know that its Fourier coefficients $\hat c_k$
satisfy $\vert \hat c_k - c_k \vert \leq \e$ for any $k$.

\medskip

3. At this step we determine the number of the equations (i.e. of
the Fourier coefficients of $f$) we need to achieve the prescribed
accuracy $\e$. We pick an appropriate finite subsystem (**) of
(*). Then we solve (**) with respect to the unknown parameters of
$\Delta P$. By the robustness estimates of \cite{Sar.Yom2} the
solution differs from the true parameters of $\Delta P$ by at most
$C_2\e$.

\medskip

4. We form a piecewise-polynomial $\hat {\Delta P}$ with the
parameters found in step 3. By the above estimates, the jump
points of $\hat {\Delta P}$ approximate the true jump points of
$f$ with the accuracy $C_3\e$ while the partial polynomials of
$\hat {\Delta P}$ approximate the values of $f$ on its continuity
intervals with the accuracy $C_4\e$. This completes the proof.

\medskip

In \cite{Sar.Yom2} we provide a detailed proof. We also compare
the above results with the classical results of Approximation
Theory on one side, and with some recent results on {\it linear
(or semi-linear)} reconstruction methods for piecewise-smooth
functions (see \cite{Tad1}-\cite{Tad4}, \cite{Vet1,Vet3,Vet5}, and
references there).

\section{Digital images}
\setcounter{equation}{0}

Our considerations in Sections 2-4 above were motivated, in
particular, by an attempt to estimate the expected efficiency of
linear versus non-linear methods methods of acquisition and
compression of still images and video-sequences.

Application of rigorous mathematical tools in Image Analysis is
usually difficult, because of an appeal to a ``human visual
perception" which is central in this field. For example, the main
compression requirement is to preserve image's ``visual quality" -
the notion which is well known to escape any attempt of a rigorous
mathematical definition.

Still, simple characteristics of images approximation, like
$L^2$-error, while not completely adequate to the ``human visual
perception" results, are usually very instructive. In the present
section we shall try to translate the rigorous results of Sections
2-4 about piecewise constant functions to the language of images.
By the reasons that become clear below we believe that our
conclusions (which we call ``statements" not to mix with theorems)
are as accurate as possible: they can be made rigorous by
restricting accurately a set of allowed images we work with.

\subsection{Linear space of images}

A typical image is represented by a rectangular array of pixels
(say, $512\times 512$). At each pixel the brightness (or the
color) discretized value is stored, typically, $8$ bits or $256$
brightness values, for grey-level images, and $24$ bits for
three-color RGB images. In this paper we shall ignore the discrete
nature of digital images, and consider them as bounded functions
on the square $Q^2$. (See \cite{Kou,Bri.Eli.Yom} for the
discussion of some specific problems related to the discrete
nature of images).

To make the space of images a linear one, we have to ignore
another important feature of true images: the image brightness has
always to be within the prescribed interval (say, [0,255]). So we
cannot add images as usual functions. Still, it is convenient to
consider images as the elements of the Hilbert space ${\mathcal
{L}}=L^2(Q^2)$ of functions on the unit two-dimensional cell
$Q^2$.

\medskip

However, considering images as elements in the linear space
$\mathcal {L}$ stresses their non-linear nature. Let us mention
some of the most immediate manifestations of this important fact.

\medskip

1. First of all, addition (or, more generally, forming linear
combinations of images) usually produces a new brightness
function, which is difficult to interpret as a meaningful
``image". Indeed, such a sum will show an artificial overlapping
of the objects appearing on each one of the original images. Only
for images representing exactly the same scene (like, for example,
the three color separations R, G, B of the same color picture)
their linear combinations have a direct visual meaning.

\medskip

2. Secondly, only a small fraction of the standard image
processing operations (like high-pass and low-pass filtering) are
linear transformations of the Hilbert space $\mathcal {L}$. Most
of the usual image processing operations (as represented, for
example, in the Adobe's ``Photoshop" package) take into account
the visual patterns on the image. Consequently, the processing is
subordinated to the geometry of the objects on the image, and in
this way it is highly non-linear.

\medskip

3. Third, individual images depend in a highly non-linear way on
the boundaries data of the objects.

\medskip

4. Finally, the most important time-dependent families of images -
video-sequences - turn out to be very complicated curves in
$\mathcal {L}$. In fact, as we shall see below, they behave
geometrically as the ``crinkled arcs" considered in Section 2.

\medskip

Let us consider in more details the effect of a motion of objects
on the image: this is the main content of typical video-sequences.
First, to simplify considerations, let us assume that the objects
are perfectly black while the background is perfectly white. Then
our images, as the elements of the Hilbert space $\mathcal {L}$,
are just the (negative) characteristic functions of the domains
occupied by the objects on the image.

If an object moves in such a way that the occupied domains are
expanding (for example, the object approaches the camera) then the
corresponding trajectory in the space of images $\mathcal {L}$ is
a crinkled arc by Proposition 2.4.

If the objects move in a more complicated way (in particular,
their boundaries are deformed in a non-rigid manner), then the
corresponding trajectory in the image space may not be a crinkled
arc. However, Proposition 2.5 shows that {\it typically} such
trajectories look like crinkled arcs ``in a small scale": the
angle between any two disjoint chords of the corresponding curve
in $\mathcal {L}$ tends to $\frac{\pi}{2}$ as these chords tend to
the same point.

The conclusion of Proposition 2.5 remains essentially valid also
under more realistic assumptions on the color of the moving
objects: indeed, near the object boundaries the image brightness
in any case behaves as a scalar multiple of the characteristic
function of the occupied domain.

Moreover, the occlusions of the moving objects do not change this
pattern. Indeed, only at the intersections of the boundaries of
the occluded objects we can expect new phenomena, but typically
these intersections have nearly zero area. So they do not affect
the $L^2$ geometry of the trajectory.

Thus we get a general (and, to our point, rather surprising)
conclusion:

\medskip

\noindent{\bf Statement 5.1} {\it A typical video-sequence is
metrically similar to a ``crinkled arc" in the Hilbert space
${\mathcal {L}}$ of images. In particular, its $\e$-entropy and
Kolmogorov $n$-width behave as those of the curve ${\mathcal
{H}}$.}

\medskip

As in Sections 2-4 above, this fact provides an immediate
limitation on the performance of linear approximation and
acquisition methods. Let us assume that we want to represent all
the images in a set $\Omega \subset {\mathcal {L}}$ which is
``large enough": namely, it contains together with each image $I$
also images representing a ``motion of objects" in $I$ (in
particular, their zoom, translations, etc.). Hence the set
$\Omega$ in fact contains ``video-sequences", and hence Statement
5.1 implies:

\medskip

\noindent{\bf Statement 5.2} {\it No $n$-dimensional linear
subspace $W$ in the Hilbert space ${\mathcal {L}}$ of images can
approximate all the images in $\Omega$ at once better than to
$C{1\over {\sqrt n}}$.}

\medskip

As for the problem of image acquisition from measurements, we get
the following conclusion, analogous to that of Theorem 4.1:

\medskip

\noindent{\bf Statement 5.3} {\it Let the image acquisition
process comprise taking $n$ measurements $m_i(I), \ i=1,\dots,n$
(linear or non-linear) of the image $I$, together with a
consequents processing $P$ of the measurements. If the processing
operator $\hat I = P(m_1,\dots,m_n)$ is a linear one then for some
$I\in \Omega$ the error $\vert \vert I - \hat I \vert \vert$ is at
least $C_1{1\over {\sqrt n}}$.}

\medskip

The inherent limitation of linear acquisition and representation
methods forced development of non-linear approaches. Most of them
utilize the fact that wavelet representations of typical images in
appropriate wavelet bases are ``sparse" - with only few ``large"
coefficients. Efficient image capturing and representation
approaches based on this fact are given, in particular in
\cite{Can.Don,Vet1,Vet3,Vet5}.

Another approach is based on ``non-linear model-net approximation"
(\cite{Bri.Eli.Yom,Eli.Yom,Yom.Zah}). In the present paper it was
described in a ``toy example" of piecewise-constant functions. As
the images are concerned, the main problem is whether such a
``model-based" representation is possible at all. Let us discuss
shortly the ``state of the art" here.

\subsection{Capturing of images and video-sequences by geometric
models}

As it was stressed in the introduction, the key to a successful
application of the ``algebraic reconstruction methods" presented
in this paper to real problems in Signal Processing lies in a
``model-based" representation of signals and especially of images.

From the point of view pursued in this paper, most of the
conventional image representation (``compression") methods can be
considered as ``semi-linear": their starting point is a linear
representation of the image in a certain basis (Fourier, local
Fourier, Wavelets ...). Then the coefficients of this linear
representation are truncated, ordered and finally encoded in a
highly non-linear way.

There are ``geometric" methods of image representation, based on
an approximation by non-linear image models (usually constructed
from the edges, ridges and other geometric visual patterns
appearing in typical images) - see
\cite{Kun.Iko.Koc,Kou,Eld,Bri.Eli.Yom,Eli.Yom}. Some of these
geometric methods have proved themselves to be very efficient in a
representation and processing of special types of images (like
geographic maps, cartoon animations, etc.).

However, in general the ``geometric" methods, as for today, suffer
from an inability to achieve a full visual quality for high
resolution photo-realistic images of the real worlds. {\it In
fact, the mere possibility of a faithful capturing such images
with geometric models presents one of important open problems in
Image Processing, sometimes called ``the vectorization problem"}.

Let us stress our strong belief that a full visual quality
``geometric" representation for high resolution photo-realistic
images of the real worlds is possible. As achieved, it promises to
bring a major advance in image compression and capturing, in
particular, via the approach of the present paper.

\subsection{Reconstruction of images from measurements}

Let a function $f(x,y)$ of two variables be the brightness
function of an image to be reconstructed from Computer Tomography
measurements. The data of the Radon transform can be translated
into the Fourier data, so we can assume that our measurements are
just Fourier coefficients of $f$.

Now our general approach to this problem extends the non-linear
inversion method presented in Section 4 above. It can be
summarized as follows:

\medskip

1. We obtain the ``first approximation" $\hat f$ of the function
$f$ to be reconstructed by one of available conventional methods.

\medskip

2. We approximate the function $\hat f$ (and hance also the
function $f$ to be reconstructed) by a ``model one" $Mf$. The last
comprises ``simple geometric models" reflecting the structure of
singularities of $f$ and approximating $f$ at its regular regions.

Let us stress once more, as we did in our discussion above, that
the mere existence of such a representation for real world images
is an important open problem.

\medskip

3. We memorize the ``combinatorial structure" of $Mf$ (the number
of its jumps in one-dimensional case; the topological structure of
the edges, ridges, patches etc. for images).

\medskip

4. We consider specific geometric and brightness parameters of the
models as unknowns, which we substitute into a system (***)
obtained in the same way as the system (4.1) above. The right-hand
side of this system is formed by the measured Fourier coefficients
of $f$.

\medskip

5. We solve the appropriate subsystem of the system (***). In the
solution process we start with the approximate solution obtained
in step 2. The solution provides a set of the improved parameters
of the model function $Mf$. Applying these parameters we finally
get an improved approximation $\hat Mf$ of the original function
$f$.

\medskip

The implementation of this program in real applications of
Computer Tomography is now in its initial stages. In ``toy
problems" where we pretend to know a priori the model structure of
the image, the approach works perfectly (not a big surprise! see
\cite{Sar.Yom1}). We believe however that the time is ripe to
study both the Image Processing and the Algebraic-Geometric parts
of the problem.

\subsection{Geometric image compression and crinkles arcs}

The following remarks concern a possibility to use directly the
universality of the ``crinkled arc" in image compression. However,
it involves encoding of certain isometries of Hilbert spaces, and
the feasibility and complexity of this task should be considered
as an absolutely open problem.

One of the most important tasks in the ``geometric image
compression" is a compact representation of the systems of curves
and points on the plane (see, in particular,
\cite{Bri.Eli.Yom,Don1,Don2}). Mostly we can assume these curves
and points to be mutually disjoint, and their specific
parametrization is not essential.

Let us consider a special case where the family of curves to be
memorized is the family of the boundaries of a family of expanding
domains in the plane. It turns to be difficult to utilize this
special structure in the curve compression methods used in
\cite{Bri.Eli.Yom, Eli.Yom}. In fact, according to these methods,
each curve will be stored {\it separately}. On the other hand, by
Propositions 2.4, 2.5 the characteristic functions of the inside
domains of our curves form a crinkled arc in $L^2(Q^2)$ isomorphic
to $\mathcal {H}$. Consequently, we have an alternative approach
to memorizing our family of curves: {\it it is enough to memorize
the transformations bringing it to $\mathcal {H}$}. This lead to
two mathematical problems which we consider as important by
themselves:

\medskip

\noindent {\bf Problem 1.} What is the complexity of the
``normalizing transformation" in Theorem 2.2, and specifically, in
Propositions 2.4, 2.5? (We can use, for example, the notion of
complexity for infinite-dimensional objects, introduced in
\cite{Yom1,Yom2}). How many bits do we need to memorize them?

\medskip

\noindent {\bf Problem 2.} How to use ``geometric redundancy" of
the expanding family - the fact that the curves do not intersect
and ``bound one another" - in their ``conventional" compression?

\bibliographystyle{amsplain}

\end{document}